\documentclass[11pt, reqno]{amsart}
\usepackage{graphics}
\usepackage{amsfonts}

\def\Box{\vcenter{\vbox{\hrule\hbox{\vrule
     \vbox to 8.8pt{\hbox to 10pt{}\vfill}\vrule}\hrule}}}

\newcommand{\GF}{{\rm GF}}

\newcommand{\Z}{{\mathbb Z}}
\newcommand{\wt}[1]{\widetilde{#1}}

\newtheorem{thm}{Theorem}[section]
\newtheorem{lemma}[thm]{Lemma}

\newtheorem{prop}[thm]{Proposition}

\DeclareMathOperator{\PG}{PG} 

\numberwithin{equation}{section}

\begin{document}

\title[Symmetric Bush-type Hadamard Matrices]
{Symmetric Bush-type Hadamard matrices of order $4m^4$ exist for all odd $m$}

\author[Mikhail Muzychuk, Qing Xiang]
{Mikhail Muzychuk, Qing Xiang$^*$}
\thanks{$^*$Research supported in part by NSF Grant DMS 0400411.}

\address{Department of Computer Sciences and Mathematics, Netanya Academic College, University St. 1, 42365, Netanya, Israel,
email: {\tt mikhail@netvision.net.il}}

\address{Department of Mathematical Sciences, University of Delaware, Newark, DE 19716, USA,
email: {\tt xiang@math.udel.edu}}


\date{}

\begin{abstract}
Using reversible Hadamard difference sets, we construct symmetric Bush-type Hadamard matrices of order $4m^4$ for
all odd integers $m$.
\end{abstract}
\maketitle

\section{Introduction}\label{intro}
A {\it Hadamard matrix} of order $n$ is an $n$ by $n$ matrix $H$ with entries $\pm 1$, such that
$$HH^{\top}=nI_n,$$
where $I_n$ is the identity matrix of order $n$. It can be easily shown that if $n$ is the order of a Hadamard
matrix, then $n=1$, 2 or $n\equiv 0$ (mod 4). The famous Hadamard matrix conjecture states that for every
positive integer $n$ divisible by 4, there exists a Hadamard matrix of order $n$. This conjecture is far from
being proved. We refer the reader to \cite{kh} for a recent construction of a Hadamard matrix of order $428$ (the
smallest order for which an example of Hadamard matrix was not known for many years). In this note, we
concentrate on a class of Hadamard matrices of highly specialized form, namely the Bush-type Hadamard matrices.

Let $n$ be a positive integer and let $J_{2n}$ denote the matrix of order $2n$ with all entries being ones. A
Hadamard matrix $H=(H_{ij})$ of order $4n^2$, where $H_{ij}$ are $2n\times 2n$ block matrices, is said to be of
{\it Bush-type} if
\begin{equation}\label{defbush}
H_{ii}=J_{2n},\; \mbox{and}\;  H_{ij}J_{2n}=J_{2n}H_{ij}=0,
\end{equation}
for $i\neq j$, $1\leq i,j\leq 2n$. K. A. Bush \cite{bush} proved that the existence of a projective plane of
order $2n$ implies the existence of a symmetric Bush-type Hadamard matrix of order $4n^2$. So if one can prove
the nonexistence of symmetric Bush-type Hadamard matrices of order $4n^2$, where $n$ is odd, then the
nonexistence of a projective plane of order $2n$, where $n$ is odd, will follow. This was Bush's original
motivation for introducing Bush-type Hadamard matrices. Wallis \cite{wa} showed that $n-1$ mutually orthogonal
Latin squares of order $2n$ lead to a symmtric Bush-type Hadamard matrices of order $4n^2$. Goldbach and Claasen
\cite{gc} also proved that certain 3-class association schemes can give rise to symmetric Bush-type Hadamard
matrices. More recently, Kharaghani and his coauthors \cite{kh00, j, jkt1, jkt2, jk} rekindled the interest in
Bush-type Hadamard matrices by showing that these matrices are very useful for constructions of symmetric designs
and strongly regular graphs. Kharaghani \cite{kh00} conjectured that Bush-type Hadamard matrices of order $4n^2$
exist for all $n$. While it is relatively easy to construct Bush-type Hadamard matrices of order $4n^2$ for all
even $n$ for which a Hadamard matrix of order $2n$ exists (see \cite{kh85}), it is not easy to decide whether
such matrices of order $4n^2$ exist if $n>1$ is an odd integer. In a recent survey \cite{jk}, Jungnickel and
Kharaghani wrote ``Bush-type Hadamard matrices of order $4n^2$, where $n$ is odd, seem pretty hard to construct.
Examples are known for $n=3$, $n=5$, and $n=9$ (see \cite{j}, \cite{jkt1}, and \cite{jkt2} respectively); all
other cases are open''. In this note, we will show that symmetric Bush-type Hadamard matrices of order $4m^4$
exist for all odd $m$.

We first note a relation between symmetric Bush-type Hadamard matrices and strongly regular graphs with certain
properties. The following lemma is well known. A weaker form of the lemma appeared in \cite{wa}. For convenience
of the reader, we provide a proof.

\begin{lemma}\label{bushSRG}
There exists a symmetric Bush-type Hadamard matrix of order $4n^2$ if and only if there exists a strongly regular
graph (SRG in short) with parameters
$$v=4n^2, \; k=2n^2-n,\; \lambda=\mu=n^2-n,$$
and with the additional property that the vertex set can be partitioned into $2n$ disjoint cocliques of size
$2n$.
\end{lemma}

\begin{proof} If $H=(H_{ij})$, where $H_{ij}$ are $2n\times 2n$ block matrices, is a symmetric Bush-type Hadamard
matrix, then the matrix $A=\frac {1}{2}(J-H)$ is symmetric and satisfies
$$A^2=n^2I+(n^2-n)J.$$
Moreover the $2n\times 2n$ block matrices on the main diagonal of $A$ are all zero matrices. Hence $A$ is the
adjacency matrix of an SRG with parameters $v=4n^2, k=2n^2-n, \lambda=\mu=n^2-n,$ and with the additional
property that the vertex set can be partitioned into $2n$ disjoint cocliques of size $2n$. Conversely, if $A$ is
the adjacency matrix of such an SRG, then the matrix $H=J-2A$ is symmetric and satisfies $H^2=4n^2I$. Since the
vertex set of the SRG can be partitioned into $2n$ cocliques, each of size $2n$, we may arrange the rows and
columns of $H$ so that we can partition $H$ into $H=(H_{ij})$, where $H_{ij}$ are $2n\times 2n$ block matrices
and $H_{ii}=J_{2n}$. It remains to show that $H_{ij}J_{2n}=J_{2n}H_{ij}=0$ for $i\neq j$, $1\leq i,j\leq 2n$.
Noting that the SRG has the smallest eigenvalue $-n$, we see that the cocliques of size $2n$ of the SRG meet the
Delsarte bound (sometime called the Hoffman bound also). By Proposition 1.3.2 \cite[p.~10]{bcn}, every vertex in
the SRG outside a coclique is adjacent to exactly $n$ vertices of the coclique. This proves that
$H_{ij}J_{2n}=J_{2n}H_{ij}=0$ for $i\neq j$, $1\leq i,j\leq 2n$. The proof is complete. \end{proof}

\section{Symmetric Bush-type Hadamard matrices from reversible Hadamard difference sets}\label{BushfromRHDS}

We start with a very brief introduction to difference sets. For a thorough treatment of difference sets, we refer
the reader to \cite[Ch.~6]{bjl}. Let $G$ be a finite group of order $v$. A $k$-element subset $D$ of $G$ is
called a $(v,k,\lambda)$ {\it difference set} in $G$ if the list of ``differences'' $d_1d_2^{-1}$, $d_1,d_2\in
D$, $d_1\neq d_2$, represents each nonidentity element in $G$ exactly $\lambda$ times. Using multiplicative
notation for the group operation, $D$ is a $(v,k,\lambda)$ difference set in $G$ if and only if it satisfies the
following equation in $\Z[G]$:
\begin{equation}\label{defdiffset}
DD^{(-1)}=(k-\lambda)1_G+\lambda G,
\end{equation}
where $D=\sum_{d\in D}d$, $D^{(-1)}=\sum_{d\in D}d^{-1}$, and $1_G$ is the identity element of $G$. A subset $D$
of $G$ is called {\it reversible} if $D^{(-1)}=D$. Note that if $D$ is a reversible difference set, then
\begin{equation}\label{defpds}
D^2=(k-\lambda)1_G+\lambda G.
\end{equation}
If furthermore we require that $1_G\not\in D$, then from (\ref{defpds}) we see that the Cayley graph ${\rm {\bf
Cay}}(G,D)$, with vertex set $G$ and two vertices $x$ and $y$ being adjacent if and only if $xy^{-1}\in D$, is an
SRG with parameters $(v,k,\lambda,\lambda)$.

The difference sets considered in this note have parameters
$$(v,k,\lambda)=(4n^2,2n^2-n,n^2-n).$$
These difference sets are called {\it Hadamard} difference sets (HDS), since their $(1,-1)$-incidence matrices
are Hadamard matrices. Alternative names used by other authors are Menon difference sets and H-sets. We will show
that reversible HDS give rise to symmetric Bush-type Hadamard matrices.

\begin{prop}\label{Hdisj}
Let $D$ be a reversible HDS in a group $G$ with parameters $(4n^2,2n^2-n,n^2-n)$. If there exists a subgroup
$H\leq G$ of order $2n$ such that $D\cap H=\emptyset$, then there exists a Bush-type symmetric Hadamard matrix of
order $4n^2$.
\end{prop}

\begin{proof} First note that the Cayley graph ${\rm {\bf Cay}}(G,D)$ is strongly regular with parameters
$v=4n^2-n,k=2n^2-n,\lambda=\mu=n^2-n.$ The cosets of $H$ in $G$ partition $G$. Let $Hg$ be an arbitrary coset of
$H$ in $G$. Then any two elements $x,y\in Hg$ are not adjacent in ${\rm {\bf Cay}}(G,D)$ since $xy^{-1}\in H$ and
$D\cap H=\emptyset$. Therefore the vertex set of ${\rm {\bf Cay}}(G,D)$ can be partitioned into $2n$ disjoint
cocliques of size $2n$. By Lemma~\ref{bushSRG}, the $(1,-1)$-adjacency matrix of ${\rm {\bf Cay}}(G,D)$ is a
symmetric Bush-type Hadamard matrix of order $4n^2$. \end{proof}

Let $G = K\times W$ where $K=\{g_0=1,g_1,g_2,g_3\}$ is a Klein four group and $W$ is a group of order $n^2$. Each
subset $D$ of $G$ has a unique decomposition into a disjoint union $D=\cup_{i=0}^3 (g_i, D_i)$ where $D_i\subset
W$. Note that if $D$ is a reversible HDS in $G$, then $(g_i,1) D$ are also reversible HDS for all $i$, $0\leq i \leq 3$. This
observation implies the following
\begin{prop}\label{Pdisj}
Let $K=\{g_0=1,g_1,g_2,g_3\}$ be a Klein four group. Let $D = \cup_{\ell=0}^3 (g_{\ell},D_{\ell})$ be a reversible Hadamard
difference set in the group $G = K\times W$, where $D_{\ell}\subseteq W$ and $|W|=n^2$. If there exists a subgroup
$P\leq W$ of order $n$ such that $P\cap D_i = P \cap D_j = \emptyset$ for some $i\neq j$, $0\leq i, j\leq 3$, then there exists a symmetric Bush-type Hadamard matrix of order $4n^2$.
\end{prop}

\begin{proof} Let $E=(g_i,1)D$. By the above observation, $E$ is a reversible HDS in $G$. Let $H=(g_0,1)P\cup
(g_ig_j,1)P$. Then $H$ is a subgroup of $G$ of order $2n$ and $H\cap E=\emptyset$. By Proposition~\ref{Hdisj},
$E$ gives rise to a symmetric Bush-type Hadamard matrix of order $4n^2$. \end{proof}

\section{Construction of symmetric Bush-type Hadamard matrices of order $4m^4$ for all odd $m$}

A symmetric Bush-type Hadamard matrix $H$ of order 4 is exhibited below.
$$H=
\begin{pmatrix}
1&1&1&-1\\
1&1&-1&1\\
1&-1&1&1\\
-1&1&1&1
\end{pmatrix}$$
So we will only be concerned with construction of symmetric Bush-type Hadamard matrices of order $4m^4$ for odd
$m>1$. We will first construct symmetric Bush-type Hadamard matrices of order $4p^4$, where $p$ is an odd prime,
from certain $(4p^4, 2p^4-p^2,p^4-p^2)$ HDS. To this end, we need to recall a construction of HDS with these
parameters from \cite{wx}. Let $p$ be an odd prime and $\PG(3,p)$ be a projective 3-space over $\GF(p)$. We will
call a set $C$ of points in $\PG(3,p)$ {\it type} Q if
$$|C|=\frac{(p^4-1)}{4(p-1)}$$
and each plane of $\PG(3,p)$ meets $C$ in either $\frac{(p-1)^2}{4}$ points or $\frac{(p+1)^2}{4}$ points. For
each set $X$ of points in $\PG(3,p)$ we denote by $\wt{X}$ the set of all non-zero vectors $v\in\GF(p)^4$ with
the property that $\langle v\rangle\in X$, where $\langle v\rangle$ is the 1-dimensional subspace of $\GF(p)^4$
generated by $v$.


Let $S=\{L_1,L_2,\ldots ,L_{p^2+1}\}$ be a spread of $\PG(3,p)$ and let $C_0, C_1$ be two sets of type Q in
$\PG(3,p)$ such that
\begin{equation}\label{eq_1}
\forall_{1\leq i\leq s}\  |C_0\cap L_i|=\frac{p+1}{2} \mbox{ and } \forall_{s+1\leq i \leq 2 s}\ \  |C_1\cap
L_i|=\frac{p+1}{2},
\end{equation}
where $s:=\frac{p^2+1}{2}$. (We note that if we take $S$ to be the regular spread in $\PG(3,p)$, then examples of
type Q sets $C_0, C_1$ in $\PG(3,p)$ satisfying (\ref{eq_1}) were first constructed in \cite{xia} when $p\equiv
3$ (mod 4), in \cite{et, wx} when $p=5, 13, 17$, and in \cite{chen} for all odd prime $p$.) As in \cite{wx} we set $C_2:=(L_1\cup ...\cup L_s)\setminus
C_0$ and $C_3:=(L_{s+1}\cup ...\cup L_{2s})\setminus C_1$. Note that $C_0\cup C_2 = L_1\cup ...\cup L_s$ and
$C_1\cup C_3 = L_{s+1}\cup ...\cup L_{2s}$.

Let $A$ (resp. $B$) be a union of $(s-1)/2$ lines from $\{L_{s+1},...,L_{2s}\}$ (resp. $\{L_1,...,L_s\}$). Let
$K=\{g_0=1,g_1,g_2,g_3\}$ and $W=(\GF(p)^4,+)$. Denote
$$
\begin{array}{rcl}
D_{0} & := & \wt{C_0} \cup \wt{A},\\
D_{2} & := & \wt{C_2} \cup \wt{A},\\
D_{1} & := & \wt{C_1} \cup \wt{B},\\
D_{3} & := & W\setminus(\wt{C_3} \cup \wt{B}).
\end{array}
$$
Then
$$|D_0|=|D_1|=|D_2|=\frac {p^4-p^2} {2},\; |D_3|=\frac {p^4+p^2} {2}.$$
By Theorem 2.2 \cite{wx} the set
$$
D := (g_0, D_{0}) \cup (g_1, D_{1}) \cup (g_2, D_{2}) \cup (g_3, D_{3}),
$$
is a reversible $(4p^4,2p^4-p^2,p^4-p^2)$ difference set in the group $K\times W$.

Pick an arbitrary line, say $L_a$, from the set $\{L_{s+1},...,L_{2s}\}$ such that $L_a\cap A=\emptyset$. Then
$P:=\wt{L_a}\cup\{0\}$ is a subgroup of $W$ of order $p^2$ such that $P\cap D_0 = P\cap D_2 = \emptyset$. Now
Proposition~\ref{Pdisj} implies that there exists a symmetric Bush-type Hadamard matrix of order $4p^4$.
Therefore we have proved

\begin{thm}\label{4p^4}
There exists a symmetric Bush-type Hadamard matrix of order $4p^4$ for every odd prime $p$.
\end{thm}

In order to build a symmetric Bush-type Hadamard matrix of order $4m^4$ for arbitrary odd $m>1$ we need to use
Turyn's composition theorem \cite{turyn}. We also need the following simple
\begin{prop}\label{Qsubgr}
There exists a subgroup $Q\leq W$ of order $p^2$ such that $Q\subseteq D_3$ and $Q\cap D_1 = \emptyset$.
\end{prop}
\begin{proof} Pick an arbitrary line $L_b$ from $\{L_1,...,L_s\}$ such that $L_b\cap B=\emptyset$ and set
$Q:=\{0\}\cup \wt{L_b}$. \end{proof}

Next we recall Turyn's composition theorem. We will use the version as stated in Theorem 6.5 \cite[p.~45]{dj}. For
convenience we introduce the following notation. Let $W_1$, $W_2$ be two groups. For $A,B\subseteq W_1$ and
$C,D\subseteq W_2$, we define the following subset of $W_1\times W_2$.
$$
\nabla({A},{B};{C},{D}):=\big((A\cap B)\times C'\big)\cup \big((A'\cap B')\times C\big)\cup \big((A\cap B')\times
D'\big)\cup\big((A'\cap B)\times D\big),
$$
where $A'=W_1\setminus A$, $B'=W_1\setminus B$, $C'=W_2\setminus C$, and $D'=W_2\setminus D$.

\begin{thm}\label{turyn}{\em (Turyn \cite{turyn})}
Let $K=\{g_0,g_1,g_2,g_3\}$ be a Klein four group. Let $E_1 = \cup_{i=0}^3 (g_i,A_i)$ and $E_2 = \cup_{i=0}^3
(g_i,B_i)$ be reversible Hadamard difference sets in groups $K\times W_1$ and $K\times W_2$, respectively, where
$|W_1|=w_1^2$ and $|W_2|=w_2^2$, $w_1$ and $w_2$ are odd, $A_i\subseteq W_1$ and $B_i\subseteq W_2$, and
$$|A_0|=|A_1|=|A_2|=\frac {w_1^2-w_1} {2},\;|A_3|=\frac {w_1^2+w_1} {2},$$
$$|B_0|=\frac {w_2^2+w_2} {2},\; |B_1|=|B_2|=|B_3|=\frac {w_2^2-w_2} {2}.$$
Define
$$
\begin{array}{rcl}
E & := & (g_0,\nabla({A_0},{A_1};{B_0},{B_1}))\cup (g_1,\nabla({A_0},{A_1};{B_2},{B_3}))\\
\ & \cup &  (g_2,\nabla({A_2},{A_3};{B_0},{B_1}))\cup (g_3,\nabla({A_2},{A_3};{B_2},{B_3})).
\end{array}
$$
Then
$$|\nabla({A_0},{A_1};{B_0},{B_1})|=\frac {w_1^2w_2^2+w_1w_2} {2},$$
$$|\nabla({A_0},{A_1};{B_2},{B_3})|=|\nabla({A_2},{A_3};{B_0},{B_1})|=|\nabla({A_2},{A_3};{B_2},{B_3})|=\frac
{w_1^2w_2^2-w_1w_2} {2},$$ and $E$ is a reversible $(4w_1^2w_2^2, 2w_1^2w_2^2-w_1w_2, w_1^2w_2^2-w_1w_2)$
Hadamard difference set in the group $ K\times W_1\times W_2$.
\end{thm}

\begin{prop}\label{intersec}
With the assumptions as in Theorem~\ref{turyn}. Let $Q\leq W_1$ and $P\leq W_2$ be such that $Q\cap A_2
=\emptyset, Q\subseteq A_3$ and $P\cap B_1 = \emptyset, P\cap B_3 =\emptyset$. Then $(Q\times P)\cap
\nabla({A_2},{A_3};{B_0},{B_1}) = \emptyset$ and $(Q\times P)\cap \nabla({A_2},{A_3};{B_2},{B_3}) = \emptyset$.
\end{prop}
\begin{proof}
It follows from the intersections
$$
\begin{array}{rcl}
Q\cap (A_2\cap A_3) & = & \emptyset,\\
Q\cap (A'_2\cap A'_3) & = & \emptyset,\\
Q\cap (A_2\cap A'_3) & = & \emptyset,\\
Q\cap (A'_2\cap A_3) & = & Q,
\end{array}
$$
that $\nabla({A_2},{A_3};{B_0},{B_1})\cap (Q\times P) = Q\times (B_1\cap P) = \emptyset$ and
$\nabla({A_2},{A_3};{B_2},{B_3})\cap (Q\times P) = Q\times (B_3\cap P) = \emptyset$.
\end{proof}

\begin{thm}\label{main}
There exists a symmetric Bush-type Hadamard matrix of order $4m^4$ for all odd $m$.
\end{thm}

\begin{proof} We only need to prove the theorem for odd $m>1$. Let $K=\{g_0,g_1,g_2,g_3\}$ be a Klein four group.
Let $p$ and $q$ be two odd primes, not necessarily distinct, and let $W_1=(\GF(p)^4,+)$ and $W_2=(\GF(q)^4,+)$.
By the construction before the statement of Theorem~\ref{4p^4}, we can construct a reversible HDS
$$E_1=(g_0, A_0)\cup (g_1, A_1)\cup (g_2, A_2)\cup (g_3, A_3)$$
in $K\times W_1$ such that
$$|A_0|=|A_1|=|A_2|=\frac {p^4-p^2} {2},\;|A_3|=\frac {p^4+p^2} {2},$$
and there exists a subgroup $Q\leq W_1$ of order $p^2$ with the property that $Q\cap A_2=\emptyset, Q\subset A_3$. (See Proposition~\ref{Qsubgr}. Note that here the $A_i$ are a renumbering of the $D_i$; any renumbering of the $D_i$ still yields a reversible difference set.) Also we can construct a reversible HDS
$$E_2=(g_0, B_0)\cup (g_1, B_1)\cup (g_2, B_2)\cup (g_3, B_3)$$
in $K\times W_2$ such that
$$|B_0|=\frac {q^4+q^2} {2},\; |B_1|=|B_2|=|B_3|=\frac {q^4-q^2} {2},$$
and there exists a subgroup $P\leq W_2$ of order $q^2$ with the property that $P\cap B_1=\emptyset, P\cap
B_3=\emptyset$ (see the paragraph before the statement of Theorem~\ref{4p^4}). Now we apply Theorem~\ref{turyn}
to $E_1$ and $E_2$ to obtain a reversible HDS
$$
\begin{array}{rcl}
E & = & (g_0,\nabla({A_0},{A_1};{B_0},{B_1}))\cup (g_1,\nabla({A_0},{A_1};{B_2},{B_3}))\\
\ & \cup &  (g_2,\nabla({A_2},{A_3};{B_0},{B_1}))\cup (g_3,\nabla({A_2},{A_3};{B_2},{B_3}))
\end{array}
$$
of size $2p^4q^4-p^2q^2$ in $K\times W_1\times W_2$. By Proposition~\ref{intersec}, we have
\begin{equation}\label{Disj}
(Q\times P)\cap \nabla({A_2},{A_3};{B_0},{B_1}) = \emptyset,\;(Q\times P)\cap \nabla({A_2},{A_3};{B_2},{B_3}) =
\emptyset.
\end{equation}
By Proposition~\ref{Pdisj}, there exists a symmetric Bush-type Hadamard matrix of order $4(pq)^4$. Now note that
$|Q\times P|=p^2q^2$, $|\nabla({A_2},{A_3};{B_0},{B_1})|=|\nabla({A_2},{A_3};{B_2},{B_3})|=\frac {p^4q^4-p^2q^2}
{2}$, and $E$ satisfies the property~(\ref{Disj}), we can repeatedly use the above process to produce a
reversible HDS satisfying the condition of Proposition~\ref{Pdisj}, hence there exists a symmetric Bush-type
Hadamard matrix of order $4m^4$ for all odd $m>1$. The proof is complete. \end{proof}

Kharaghani \cite{kh00, kh01} showed how to use Bush-type Hadamard matrices to simplify Ionin's method \cite{io}
for constructing symmetric designs. Based on his constructions in \cite{kh00, kh01}, we draw the following
consequences of Theorem~\ref{main}.

\begin{thm}
Let $m$ be an odd integer. If $q=(2m^2-1)^2$ is a prime power, then there exists twin symmetric designs with
parameters
$$v=4m^4\frac {(q^{\ell +1}-1)} {q-1},\;k=q^{\ell}(2m^4-m^2),\;\lambda=q^{\ell}(m^4-m^2),$$
for every positive integer $\ell$.
\end{thm}

\begin{thm}
Let $m$ be an odd integer. If $q=(2m^2+1)^2$ is a prime power, then there exists Siamese twin symmetric designs
with parameters
$$v=4m^4\frac {(q^{\ell +1}-1)} {q-1},\;k=q^{\ell}(2m^4+m^2),\;\lambda=q^{\ell}(m^4+m^2),$$
for every positive integer $\ell$.
\end{thm}

\end{document}